\documentclass[11pt]{amsproc}
\usepackage{fancyhdr}
\usepackage{latexsym,amsfonts,amssymb,amsthm,amsmath}
\usepackage[alphabetic]{amsrefs}
\usepackage[top=1.35in,bottom=1.35in,left=1.6in,right=1.6in]{geometry}
\usepackage[all]{xy}
\usepackage{todonotes}
\usepackage{tikz-cd}
\usepackage{tikz}
\usepackage{enumerate}

\fancypagestyle{firststyle}{
	
	\fancyhf{}
	\fancyhead[C]{\footnotesize \textit{Enchiridion: Mathematical User's Guides} \hfill Volume 3 (September 2017) \\ \url{mathusersguides.com} \hfill \copyright 2017 CC BY-NC 4.0}
}

\newtheorem{lemma}{Lemma}[section]

\newtheorem{idea}[lemma]{Key Idea}  

\theoremstyle{definition}

\theoremstyle{plain}
\newtheorem{thm}{Theorem}
\newtheorem*{thm*}{Theorem}

\newcommand{\ints}{\mathbb{Z}}
\newcommand{\reals}{\mathbb{R}}

\newcommand{\sph}{\mathbb{S}}
\newcommand{\E}{\mathbb{E}}

\newcommand{\C}{\mathcal{C}}

\title[User's guide: Relative Thom spectra]{A user's guide: Relative Thom spectra via operadic Kan extensions}

\author{Jonathan Beardsley}
\address{Department of Mathematics \\ University of Washington
\\ Seattle, WA 98195}
\email{jbeards1@uw.edu}

\begin{document}

\maketitle
\thispagestyle{firststyle}

\tableofcontents


\section{Key insights and central organizing principles} \label{section-ideas}

The main idea of \cite{beardsrelthom} is to give a ``Third Isomorphism Theorem" for quotients of ring spectra by actions of $n$-fold loop spaces. Recall the classical Third Isomorphism Theorem:

\begin{thm}[Noether]
Let $N\subseteq K\subseteq G$ be a composition of inclusions of normal subgroups. Then there is an isomorphism of groups: $$G/K\cong(G/N)/(K/N).$$ 
\end{thm}

We can generalize this theorem to any action of a Lie group on a (real or complex) smooth manifold:

\begin{thm}[Bourbaki]
Let $G$ be a Lie group acting freely, properly and smoothly on a finite dimensional smooth $X$. Let $H$ be a normal Lie subgroup of $G$. Then the canonical projection map $X\to X/H$ defines, upon taking quotients by $G$, an isomorphism of smooth manifolds: $$X/G\cong (X/H)/(G/H).$$
\end{thm}

In general, we see that when one has an action of a group $G$ on an object $X$, and a normal subgroup $H\subset G$, it is reasonable to expect that $X/G$ and $(X/H)/(G/H)$ will be very similar to each other. 

In \cite{abghr} it was shown that Thom spectra (e.g.~the classical cobordism spectra) can be described as quotients of ring spectra by the actions of $n$-fold loop spaces. Thus the question arises: given an $n$-fold loop space $G$ acting on a ring spectrum $R$, and ``normal subgroup" $H\subseteq G$, what is the relationship between $R/G$ and $(R/H)/(G/H)$? Of course to even answer this question one must make rigorous a number of notions. Namely, subgroups of $n$-fold loop spaces, actions of $n$-fold loop spaces on ring spectra and quotients of ring spectra by such actions. It's also going to be important that $n$-fold loop spaces are algebras for the little $n$-cubes operad, and as such admit monoid structures up to coherent homotopy (where the level of coherence increases with $n$).

\subsection{``Normal subgroups" of $n$-fold loop spaces}

Given a discrete group $G$ with a subgroup $H$, the most important thing one gets from normality of $H$ is the ability to equip the quotient set $G/H$ with a group structure. One might even say that this is the defining property of being normal; this is the property that makes normality worth knowing about. In the case of an $n$-fold loop space, we might similarly ask for $H$ to have some property such that $G/H$ is still an $n$-fold loop space. Completely characterizing such sub-objects is very hard, but we are rescued by a preponderance of natural examples of this structure. 

Our analogy for the inclusion of a normal subgroup $H\subseteq G$ will be a fiber sequence of $n$-fold loop spaces $F\to E \to B$, where the maps are maps \textit{of} $n$-fold loop spaces. In this analogy, $F$ will play the role of $H$, $E$ the role of $G$, and $B$ the role of $G/H$. Recall that for a morphism of discrete groups $\phi\colon G\to K$ we can take the kernel, the subgroup of $G$ that goes to 1 under $\phi$. Categorically, this is equivalent to taking the pullback of the cospan $0\to G\leftarrow K$. In the case that $K$ is $G/H$ for some normal subgroup, we have that the kernel of $\phi$, the pullback of that cospan, is isomorphic to $H$. Thus replacing $G$ and $K$ with $n$-fold loop spaces, and replacing that pullback with a homotopically coherent pullback (e.g. taking the pullback in the quasicategory of $n$-fold loop spaces), we generalize the notion of normal subgroup to $n$-fold loop spaces.

\begin{idea}
When working with $n$-fold loop spaces instead of strict groups, we should replace the notion of ``$H$ is a normal subgroup of $G$ and $G/H\cong K$" with ``there is a homotopy fiber sequence of $n$-fold loop spaces $H\to G\to K$."
\end{idea}

Probably the most obvious examples of such structure are just the loop space functor applied $n$ times to any fiber sequence of simplicial sets. For instance, we might take the well known Hopf-fibration $S^1\to S^3\to S^2$, and take some number of based loops on that fibration. Thus we would have that $\Omega^n S^1$ is a ``normal sub-$n$-fold loop space" of $\Omega^n S^3$, and $\Omega^nS^3/\Omega^nS^1=\Omega^nS^2$. Recall that we also have a number of compact Lie groups that are infinite loop spaces, and there are many infinite loop maps going between them. Thus we get fiber sequences of infinite loop spaces like $U\to O\to O/U$, $SU\to U\to S^1$. We will see several other such fiber sequences in later sections when we discuss examples. 

\subsection{Actions of $n$-fold loop spaces on ring spectra}

To discuss actions of $n$-fold loop spaces on ring spectra, we need to decide on a notion of what we mean by ring spectra. Because it was the notion used in \cite{beardsrelthom}, we'll work with the symmetric monoidal quasicategory of ring spectra described in \cite{ha}. Ultimately the choice of model of ring spectra seems inessential, but there might be difficulties that this author isn't aware of.  We will also assume that our ring spectra admit multiplicative structures parameterized by \textit{operads} (more precisely, the $\infty$-operads of \cite{ha}), in particular, the little $n$-disk operads, $\E_n$, of Boardman and Vogt \cite{boardvogt}. It would be a distraction to describe operads in detail here, but suffice it to say that $\E_1$ things should be thought of as \textit{homotopy associative with all higher coherence data}, $\E_\infty$ should be thought of as \textit{homotopy associative and commutative with all higher coherence data}, and $\E_n$ should be thought of as \textit{homotopy associative with all higher coherence data and homotopy commutative with some higher coherence data}, where ``some" will increase as $n$ increases. There is a close connection between $n$-fold loop spaces and $\E_n$-ring spectra: in particular the suspension spectrum of an $n$-fold loop space is always an $\E_n$-ring spectrum.

Generally, to define an action of a monoid $G$ on an object $X$, one provides a map of monoids: $f:G\to Aut(X)$, and that's exactly what we'll do here. Specifically, we know from \cite{abg} (as well as plenty of other places) that if $R$ is an $\E_n$-ring spectrum then there is an $n$-fold loop space of homotopy automorphisms of $R$ as an $\E_n$-ring, which we'll denote by $GL_1(R)$. So we'll say that an $n$-fold loop space $G$ acts on an $\E_n$-ring spectrum $R$ if there is a morphism of $n$-fold loop spaces $f\colon G\to GL_1(R)$. We will often want to go between the map $f\colon G\to GL_1(R)$ and $Bf\colon BG\to BGL_1(R)$ by looping and delooping, so we will always assume that our $n$-fold loop spaces are connected. 

\subsection{Quotients of ring spectra by actions of $n$-fold loop spaces}

Let's suppose we have an $\E_n$-ring spectrum $R$ and a map of $(n-1)$-fold loop spaces $f\colon BG\to BGL_1(R)$, hence an action of an $n$-fold loop space $G$ on $R$. By thinking of left $R$-modules as a quasicategory $LMod_R$, we can think of $BGL_1(R)$ as a sub-quasicategory of $Mod_R$ in the following way: the base point of $BGL_1(R)$ is $R$ itself, paths in $BGL_1(R)$ are the obvious homotopy automorphisms of $R$, higher cells in $BGL_1(R)$ are mapped to the necessary higher degree morphisms in $LMod_R$. In other words there is an inclusion of simplicial sets $BGL_1(R)\hookrightarrow LMod_R$ which precisely picks out the action of $GL_1(R)$ on $R$ as a left $R$-module. Thus a morphism of $(n-1)$-fold loop spaces $Bf\colon BG\to BGL_1(R)$ can be thought of as picking out a $BG$-shaped diagram inside of $LMod_R$ describing an action of $G$ on $R$ by left $R$-module automorphisms. 

To construct the ``quotient" of $R$ by this action in a homotopically invariant way, we should imagine forcing all of the morphisms in the image of $BG$ to be equivalent to the identity morphism. One might imagine manually attaching 2-cells between each morphism and the identity morphism, and then attempting to make this coherent. Unfortunately this would require giving an infinite list of coherences, and coherences of coherences, and coherences of coherences of... ad infinitum. Luckily, the theory of \cite{htt} and \cite{ha} builds all coherence data into the category theory itself, making this completely formal: we will say that the quotient of $R$ by the $G$ action is the \textit{quasicategorical} colimit (i.e. a quasicategorical analog of the classical notion of a homotopy colimit, hence homotopy invariant), in $LMod_R$, of the morphism $BG\to BGL_1(R)\hookrightarrow LMod_R$. Our intuition about classical categorical colimits is the right one in this case and $R/G$, or $colim(BG\to BGL_1(R)\hookrightarrow LMod_R)$,  is precisely the universal left $R$-module on which $G$ acts homotopically trivially. 

\begin{idea}
An action of an $n$-fold loop space $G$ on an $\E_n$-ring spectrum $R$ is the data of a functor $f:BG\to LMod_R$ that takes the unique basepoint of $BG$ to $R$ and every path in $BG$ to an $\E_n$-algebra automorphism of $R$. The quotient of $R$ by such an action is the $R$-module spectrum $colimit(f)$. 
\end{idea}

Now that we have all of the relevant notions in place, we can state the main theorem of \cite{beardsrelthom}:

\begin{thm*}
Suppose $Y\overset{i}\to X\overset{q}\to B$ is a fiber sequence of reduced $\E_n$-monoidal Kan complexes for $n>1$ with $i$ and $q$ both maps of $\E_n$-algebras. Let $f\colon X\to BGL_1(\sph)$ be a morphism of $\E_n$-monoidal Kan complexes for $n>1$. Then there is a a morphism of $\E_{n-1}$-algebras $B\to BGL_1(M(f\circ i))$ whose associated Thom spectrum is equivalent to $Mf$. 
\end{thm*}

Note that in the actual statement we have to be more precise then we've been up to this point. In particular, instead of talking about $n$-fold loop spaces, we specifically talk about Kan complexes which are algebras for the $\E_n$-operad. We also assume that our Kan complexes are reduced. This means that as simplicial sets they have a single zero simplex. This guarantees, for instance, that they are connected, but it also makes things technically simpler. In the end, any time one has a connected Kan complex, one can replace it up to homotopy with a reduced one. Note also that in the statement of the theorem we write $Mf$ and $M(f\circ i)$ instead of $\sph/\Omega Y$ or $\sph/\Omega X$. This is purely a notational convention and is done stay aligned with classical descriptions of such quotient spectra as Thom spectra.

\subsection{Important examples of quotients and iterated quotients}

It turns out that there are many examples in homotopy theory of taking quotients of ring spectra by $n$-fold loop space actions. As described above, these are typically called Thom spectra. Recall that there is a map of infinite loop spaces $j\colon BO\to BGL_1(\mathbb{S})$ called the $j$-homomorphism. So any time we have a map of infinite loop spaces (or just of $n$-fold loop spaces, since we can just forget some structure) $BG\to BO$, we can take the induced quotient $\mathbb{S}/G$. Thus we can obtain all of the classical Thom spectra in this way: $MO$, $MSO$, $MU$, $MSU$, $MSp$, $MSpin$ and so forth. However, it is also known that one can produce the Eilenberg-MacLane spectra $H\ints$, $H\ints/p$ and $H\widehat{\ints}_p$ in this way (though sometimes these require taking quotients of the $p$-local sphere $\sph_{(p)}$ rather than just $\sph$). These last three examples are originally due to Mahowald and Hopkins but a modern treatment can be found in \cite{blumthh}. In particular, there is map $f:\Omega^2S^3\to BGL_1(\sph_{(p)})$ whose associated quotient is $H\ints/p$ and there is a map $g:\Omega^2S^3\langle 3\rangle\to BGL_1(\sph_{(p)})$ which is factored by $f$ and has associated quotient $H\widehat{\ints}_p$. 
 
For \emph{iterated} quotients, we can first notice that all of the classical cobordism spectra described above come from infinite loop spaces maps $BG\to BO$ and that most of these factor through $BSO$. So we can think of $BG$ as being the fiber of a map $BSO\to B(SO/G)$. Thus we get that $\sph/O\simeq (\sph/G)/(SO/G)$. For technical reasons it will help to have all of our groups connected, so we'll always work with $SO$ instead of $O$. There are ways to get around this problem, but it wouldn't be productive to spend time on them in this user's guide. As an example, recall that there is a fiber sequence of infinite loop spaces $BSpin\to BSO\to B\reals P^\infty\simeq K(\ints/2, 2)$. Thus we have that $\sph/SO=MSO\simeq (\sph/Spin)/(SO/Spin)\simeq MSpin/\reals P^\infty$. 

We get another interesting example by considering the fiber sequence of 2-fold loop spaces $\Omega^2S^3\langle 3\rangle\to \Omega^2S^3\to \mathbb{C}P^\infty\simeq S^1$ along with the map $f:\Omega^2S^4\to  BGL_1(\sph_{(p)}$ described in \cite{blumthh}. In this way, we see that $H\ints/p\simeq H\widehat{\ints}_p/\Sigma^\infty_+\ints$ where the $\Sigma^\infty_+\ints$ action is of course given by multiplication by $p$. One especially nice benefit of this description is that the above construction is obtained as an equivalence of $\E_1$-ring spectra. 

Finally, near and dear to the heart of this author, is the example of the sequence of spectra $X(n)$ defined by Ravenel in \cite{rav} and used to prove the  Nilpotence Theorem in \cite{devhopsmith}. Each $X(n)$ is obtained from a map of 2-fold loop spaces $\Omega SU(n)\to BGL_1(\sph)$, i.e. $X(n)=\sph/\Omega^2SU(n)$. Recalling that there are fiber sequences $\Omega SU(n)\to \Omega SU(n+1)\to \Omega S^{2n+1}$, we have that $X(n+1)\simeq X(n)/\Omega^2S^{2n+1}$. 

\subsection{So what?}

One really interesting consequence of all this comes from the fact that whenever we have a $k$-fold loop map $BG\to BGL_1(R)$ for $R$ an $\E_{n}$ ring spectrum (with $n>k$), we get a so-called Thom isomorphism of spectra $R/G\wedge_{R} R/G\simeq R\wedge_\sph \Sigma^\infty BG_+$. So for the examples above we obtain equivalences of spectra like: $H\ints/p\wedge_{H\widehat{\ints}_p}H\ints/p\simeq H\ints/p\wedge_{\sph} S^1$ and $MSpin\wedge_{MSO} MSpin\simeq MSpin\wedge_{\sph} \Sigma^\infty_+ K(\ints/2,2).$ Thus the methods described in \cite{beardsrelthom} give relative K\"unneth theorems for many well known Thom spectra. 

Additionally, the main theorem above can be thought of as one direction of a generalized Galois correspondence. In particular, recall that in a Galois extension $E\to F$, intermediate Galois extensions are in bijection with normal subgroups of $Gal(E/F)$. In our case, we're showing that if we have a fiber sequence of $n$-fold loop spaces $H\to G\to G/H$ and a map of $(n-1)$-fold loop spaces $BG\to BGL_1(R)$, then we get a composition of spectra $R\to R/H\to R/G$. It is not true that we can recover $R/H$ from $R/G$ or $R$ from $R/H$ by taking homotopy fixed points, but we can often perform these recoveries by taking so-called homotopy \emph{co}fixed points. Thus this composition is not an iterated Galois extension but an iterated Hopf-Galois extension in the sense of \cite{rog}. So our main theorem gives a distinctly algebro-geometric interpretation to many classical (and geometric!) morphisms of cobordism spectra. 


\section{Metaphors and imagery} \label{section-metaphors}

Presenting Thom spectra as quotients of ring spectra by actions of $\E_n$-spaces, as in \cite{abghr}, ends up meaning that the most important metaphors for us involve thinking about groups acting on things. I'll explain below how I think about groups (or $\E_n$-spaces) acting on spectra, how I think about taking quotients by these actions, and finally how I think about quasicategories in general.

\subsection{Group actions as functors}

Given an $\E_n$-monoidal space $G$, we can form its classifying spaces $BG$, which is an $\E_{n-1}$-monoid and has the property that $\Omega BG\simeq G$. There are a number of ways to construct $BG$, including the well-known bar construction, but these are not always intuitively enlightening. However, since $\Omega BG\simeq G$, we at least have a good way to think about $\pi_1(BG)$. That is to say, a \emph{path} in $BG$ that starts and ends at the base point corresponds to multiplication by a point in $G$. Given a group (or $\E_n$-space) $G$ with points $\{1,\alpha,\beta, \gamma\}$, I think of $BG$ as something like the following picture: 

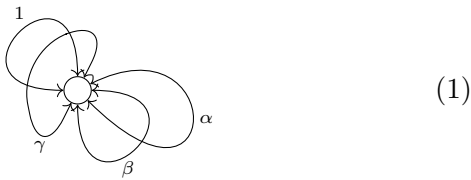
\begin{center}
\hspace{1in}
\begin{tikzpicture}
\node (group) at (5,0) [circle, draw]{};
\draw [<->] (group) .. controls (7,1) and (7,-2) ..  node[xshift=5pt]{\tiny $\alpha$} (group);
\draw [<->] (group) .. controls (7,0) and (5,-2) ..  node[below=1pt, xshift=-3pt]{\tiny $\beta$} (group);
\draw [<->] (group) .. controls (3,0) and (5,2) ..  node[above=1pt]{\tiny 1} (group);
\draw [<->, rounded corners=20pt] (group) -- (5.5,1) -- (4.2,.5) -- (4.5,-1) node[yshift=7pt]{\tiny $\gamma$} -- (group);
\node (label) at (10,0) {(1)};
\end{tikzpicture}
\end{center}
In other words, $BG$ looks like a single connected component with paths attached to it for each point of $G$. Note that $BG$ is a space (or simplicial set) and has higher homotopy groups, so the above picture does not show all of the structure of $BG$. For example, if $G$ had a relation like $\alpha\gamma^{-1}=\beta$, one would have to glue in a 2-cell to $BG$ to manifest this relation.

If we think of $BG$ as a category (or a quasicategory) it would have a single object and a 1-morphism for each point of $G$. Importantly, since we can write $G\simeq \Omega BG$, we know that every point in $G$ admits an ``inverse" corresponding to traversing the loop in the opposite direction. Thus since multiplication by a point in $G$ is always invertible (up to homotopy), $BG$, thought of as a category, is actually a groupoid (or $\infty$-groupoid). 

Suppose now that we have an $\E_{n-1}$-monoidal quasicategory $\C$ and an $\E_{n-1}$-monoidal functor $f:BG\to \C$. Since $BG$ only has one object $\ast$, the functor $f$ must pick out a diagram in $\C$ that looks like a single object and a collection of equivalences of that object (as well as higher coherent cells). Moreover, since the functor $f$ is $\E_{n-1}$-monoidal it must take $\ast$ to an $\E_n$-algebra in $\C$ and each loop in $BG$ to a morphism of $\E_{n-1}$-algebras. It also has to take 1 to the identity morphism (or at least a morphism equivalent to the identity), which will be relevant in understanding the colimit of this diagram.  Graphically, if $f(\ast)=X$,  then the image of $f$ looks like a diagram in $\C$ of the following form:

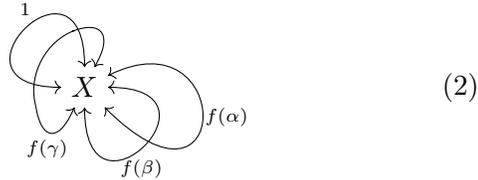
\begin{center}
\hspace{1in}
\begin{tikzpicture}
\node (group) at (5,0) {$X$};
\draw [<->] (group) .. controls (7,1) and (7,-2) ..  node[xshift=9pt]{\tiny $f(\alpha)$} (group);
\draw [<->] (group) .. controls (7,0) and (5,-2) ..  node[below=1pt, xshift=-1pt]{\tiny $f(\beta)$} (group);
\draw [<->] (group) .. controls (3,0) and (5,2) ..  node[above=1pt]{\tiny 1} (group);
\draw [<->, rounded corners=20pt] (group) -- (5.5,1) -- (4.2,.5) -- (4.5,-1) node[yshift=5pt]{\tiny $f(\gamma)$} -- (group);
\node (label) at (10,0) {(2)};
\end{tikzpicture}
\end{center} 

So the functor $f:BG\to \C$  can be though of as picking out an \emph{action} of $G$ on $X$. Given a point in $G$, $f$ tells us exactly how to use it to transform $X$ by choosing an equivalence of $X$ in $\C$ associated to it. What's also really great about this formalism is that by asking for $f$ to be a functor of quasicategories, it retains all the (potentially very complicated) coherence data contained in $G$.  

\subsection{Quotients of actions as colimits}

Now we can elaborate on what it means to take the ``quotient" of a $G$-action when the $G$ action is described as above. There are two competing pictures here that I use, depending on how I'm thinking about this quotient. We touched on both of them in Topic 1, but we'll flesh them out further here.

\subsubsection{A categorical description}
The first one is the simplest (and the most ``correct"), but also perhaps less intuitive. If we return to picture (2) above, and remember the idea behind the colimit of a diagram in $\C$, we know that $colim(f(BG))$ should be a single object of $C$ that admits a morphism \emph{from} the diagram $f(BG)$. So we have a picture like this: 

\begin{center}
\hspace{1in}
\begin{tikzpicture}
\node (group) at (2,0) {$X$};\node (quotient) at (7,0) {$X/G$};\node (secret) at (4,0) {};
\draw [<->] (group) .. controls (4,1) and (4,-2) ..  node[xshift=9pt]{\tiny $f(\alpha)$} (group);
\draw [<->] (group) .. controls (4,0) and (2,-2) ..  node[below=1pt, xshift=-1pt]{\tiny $f(\beta)$} (group);
\draw [<->] (group) .. controls (0,0) and (2,2) ..  node[above=1pt]{\tiny 1} (group);
\draw [<->, rounded corners=20pt] (group) -- (2.5,1) -- (1.2,.5) -- (1.5,-1) node[yshift=5pt]{\tiny $f(\gamma)$} -- (group);
\draw [->] (group) .. controls (3,1) and (5,1) .. node[yshift=4pt]{\tiny $q$}(quotient);
\node (label) at (10,0) {(3)};
\end{tikzpicture}
\end{center}

But the really important thing about the fact that $X/G$ is the colimit of this diagram is that this new diagram, shown in (3) above, has to commute. In other words, all possible concatenations of loops in  picture (2) followed by the morphism $q\colon X\to X/G$ have to be equivalent. And moreover, any \emph{other} object admitting a morphism from picture (2) that commutes with all those loops has to admit a map from $X/G$ factoring that morphism. In particular, for any two loops $f(\alpha)$ and $f(\beta)$ in $f(BG)$ we have to have that $q\circ f(\alpha)\simeq q\circ f(\beta)$. So, in particular, $q\circ f(\alpha)\simeq q\circ 1$. Thus $X/G$ really is the universal thing admitting a map from $X$ on which the $G$-action is forced to be homotopically trivial. 

\subsubsection{A topological description}

The second way of thinking about this quotient is a bit more ``hand-wavy." Since all of this is happening in the setting of quasicategories, which are a special type of simplicial set, we can use some of our intuition about working with simplicial sets. In particular, we can think about \emph{attaching cells} to make the loops in $f(BG)$ homotopically trivial. 

If we were working with a \emph{set} $X$ with an action by a group $G$, we would take the orbits of $G$-action on $X$, or the quotient of $X$ by $G$, by describing the quotient as being the things in $X$ after $x$ has been forced to be \emph{equal} to $gx$ for every $x\in X$ and $g\in G$.  When doing homotopy theory, and especially higher category theory, it's considered \emph{evil} to ask for things to actually be equal. Instead, if we want two things to be ``the same" in some way, we put a cell between them. This is why, for instance, if we have a based space $Y$ and take the pushout of $\ast\leftarrow Y\to \ast$, we get $\ast$, but if we take the \emph{homotopy pushout} we get $\Sigma Y$. The reduced suspension of $Y$ is precisely two points with every way of identifying them being given its own 1-cell $\ast\leftrightarrow\ast$.

For our case, we have autoequivalences $f(\alpha)\colon X\to X$ that we want to trivialize. One way to do this would be to ``glue" cells to $X$ for each $x$ and $gx$. But since we're working in quasicategories (i.e. $\infty$-categories), we would then have to glue in 2-cells between certain 1-cells, and $n$-cells in general, according to the structure of $G$. Thus while thinking of the procedure as being ``gluing in cells" we're forced to take the abstract colimit described above.

\subsection{Twisting and untwisting}

Despite the fact that we have focused so far on the idea of taking quotients of group actions, there is also a strong undercurrent of bundle theory, and twisted bundles, throughout \cite{beardsrelthom}. Recall that for any group $G$ and a space $X$, a map $X\to BG$ always determines a principal $G$-bundle over $X$. Thus in our case, a map $BG\to BGL_1(\sph)$ defines a principal $GL_1(\sph)$-bundle on $BG$. Note that if we replaced $\sph$ with a commutative ring $R$, we'd equivalently have a bundle of $R$-modules over $BG$ and this bundle would locally be isomorphic to $R$. Since it doesn't make sense to talk about a bundle of \emph{spectra} on $BG$ (for instance, what would the total space be?), we can't say the same here, but we still \textit{think} of a map $BG\to BGL_1(\sph)$ as defining a ``bundle of 1-dimensional $\sph$-modules over $BG$." 

Suppose now that the map of interest is the trivial one, $\ast\colon BG\to BGL_1(\sph)$. Then the induced quotient, $\sph/G$ where $G$ acts trivially, is in fact equivalent to $\sph\wedge\Sigma^\infty_+BG=\sph[BG]$, the spherical group ring on $BG$. Compare this with the classical fact that if $X$ is a $G$-space with a trivial $G$-action then the Borel construction $X\times_G EG\simeq X\times BG$. But now notice that if, on the other hand, we pulled back a ``bundle of sphere spectra" along the trivial map $\ast\colon BG\to BGL_1(\sph)$, we should get something that \emph{looks like} $BG\times \sph$. Of course that last expression doesn't make sense since $BG$ and $\sph$ aren't even in the same category, but this is still a useful perspective. 

For a non-trivial map, we can still use this perspective, but we should think of the map $f\colon BG\to BGL_1(\sph)$ as \emph{twisting} the trivial bundle on $BG$. And then we should think of the colimit of $ BG\overset{f}\to BGL_1(\sph)\to Spectra$ as being a \emph{twisted tensor product} of $BG$ with $\sph$. This construction should be thought of as being analogous to the twisted tensor products that arise when one attempts to compute the homology groups of fibered spaces, as in \cite{brown}. 

\subsection{Quasicategories}

Throughout \cite{beardsrelthom} and in many of the references given therein, standard category theory has been replaced by the theory of quasicategories, a model for so-called $\infty$-categories. Though there isn't space to truly explicate what makes the theory of quasicategories go, I thought it might help to at least provide a picture here of what's going on, and why one might choose to use quasicategories instead of model categories. 

Central to the entire notion of category theory is the idea that there are morphisms between objects and that one can compose a morphism $f$ with a morphism $g$ if the domain of $g$ is equal to the domain of $f$. When we do that, we get a commutative diagram like the following: 

\begin{center}
\begin{tikzcd}
&Y\arrow[dr,"g"]&\\
X\arrow[ur,"f"]\arrow[rr,"g\circ f"]&&Z
\end{tikzcd}
\end{center}

Though it seems obvious, notice that $f\circ g$ is a \emph{new} morphism whose existence we've demanded (otherwise we don't have a category). Moreover, we've demanded that the application of $g\circ f$ to $X$ is \emph{equal} to the application of $f$ followed by the application of $g$. 

For homotopy theorists, asking for something to be equal is a real problem. The main issue is that we might have two categories that are Quillen equivalent (hence have the same ``homotopy theory") and things that are equal in one might be only \emph{equivalent} in the other. Thus is general it's too restrictive to demand that compositions be equal. A good example to have in mind here is when $f$ and $g$ are continuous morphisms of spaces. We could ask that $g\circ f(x)=g(f(x))$ for every $x\in X$, or we could simply ask that there is a continuous family of continuous maps $H(t)$ such that $H(0)(x)=g\circ f(x)$ and $H(1)(x)=g(f(x))$. 

The way we get around this problem, which is intrinsic to category theory itself, is to replace categories altogether with simplicial sets. Just like categories, simplicial sets still have objects (0-simplices) and morphisms (1-simplices), but they also have a lot of ``higher data" that we think of as $n$-morphisms for $n>1$. Notice that in an arbitrary simplicial set, it really doesn't make sense to ``compose" morphisms. The 1-cells in a simplicial set $X$ are just elements in some set, specifically they're elements of $X([1])$. So to make simplicial sets behave like categories we instead do the following: given two 1-simplices $f$ and $g$ such that the right endpoint of $f$ is the same 0-simplex as the left endpoint of $g$, we can ask that there be another 1-simplex $h$ going from the left endpoint of $f$ to the right endpoint of $g$. In other words, we can ask for a triangle in our simplicial set: 

\begin{center}
\begin{tikzcd}
&y\arrow[dr,"g"]&\\
x\arrow[ur,"f"]\arrow[rr,"h"]&&z
\end{tikzcd} 
\end{center}

So we have that $h$ is modeling the \emph{composition} of $g$ with $f$. We need one more condition however: we ask that there is a 2-simplex whose boundary is the above triangle. What this is supposed to be modeling is that going along $h$ is ``the same" as going along $f$ and then along $g$. But it's not actually \emph{equal.} That two cell is taking the place of the family of continuous maps from our example. 

Visually, I think of the 2-simplex described above as giving me a way of sliding $h$ up to the edge of the triangle made by $f$ and $g$. And the fact that I'm able to do this is telling me that while $h$ may not be equal to $g(f)$, we at least have a coherent way of wiggling the data of $h$ to the data $f$ and $g$. If I can find an $h$ and a filling 2-simplex every time I have two 1-simplices that match up the way $f$ and $g$ do, then my simplicial set is behaving a lot like a category. I can always produce a ``composition" of morphisms, and this composition can always at least be homotoped back to the things whose composition it represents.

Finally, note that simplicial sets have $n$-simplices for all $n$. So asking for a simplicial set to actually be a quasicategory is asking that we can perform that above procedure in all dimensions. In other words, whenever I have 3 composable 2-simplices (think of these as forming a tetrahedron with no bottom side), I can find another 2-simplex that closes the tetrahedron and then, crucially, I can fill in that whole tetrahedron with a 3-simplex. If we can do this \emph{every time} we have $n+1$ composable $n$-simplices then our simplicial set deserves to be called a quasicategory or $\infty$-category.


\section{Story of the development} \label{section-story}

Before we dive in to the story of this particular paper, it's necessary to give a bit of background. When I got to graduate school in 2010, I only knew that I liked category theory and that Grothendieck seemed cool. After my first year I asked Jack Morava to be my advisor, primarily because, like Grothendieck, he seemed cool. What I didn't know at the time was that Jack, with his visionary ideas about the structure of the stable homotopy category, was effectively the father of what came to be called \emph{chromatic homotopy theory}. The original motivation for my thesis, and thus \cite{beardsrelthom}, came entirely from this chromatic realm.

Chromatic homotopy theory amounts to the realization that the stable homotopy groups of spheres, and subsequently the category of finite CW complexes, can be \emph{stratified} into layers that look like the height stratification on the moduli stack of formal groups. These strata are the ``colors" that give chromatic homotopy theory its name. Ravenel had made several conjectures about this structure in \cite{rav3}, the most famous being the so-called Nilpotence, Periodicity and Thick Subcategory conjectures. These were later proven in \cite{devhopsmith} and \cite{hopsmith}. 

The only thing the casual reader needs to know about these conjectures, which are now theorems, is that they: (i) tell us that if $X$ is a finite cell complex and $\alpha\in\pi_\ast(X)$ then $\alpha$ is nilpotent if and only if $MU_\ast(\alpha)$ is trivial, and (ii) the homotopy category of finite cell complexes is stratified into so-called ``thick" subcategories that directly mirror the structure of the moduli stack of formal group laws. It's this last statement that always struck me as somewhat miraculous. We have two things that are by all appearances completely unrelated: finite cell complexes and group structures on formal schemes (which come up in number theory). The Thick Subcategory Theorem however tells us that there is a deep (and still mysterious!) link between them. Understanding this connection is the true motivation for \cite{beardsrelthom}. 

As for the paper that I wrote, its genesis depends on the fact that I had a thesis advisor that more or less let me do whatever I wanted. This was both a blessing and a curse. It allowed me to avoid doing any of the really unpleasant homological algebra that is the bread and butter of algebraic topology but it also meant that after about 4 years in graduate school and 3 failed thesis projects, I was in a bit of a bind (one reason being that I wasn't very good at hard homological algebra...). I knew that I wanted to try to say something about stable homotopy theory and its relationship to algebraic geometry via the Nilpotence Theorem of \cite{devhopsmith}, but I didn't know how.

Moderately frantic conversations with my advisor and my de facto co-advisor Andrew Salch, led me to the following question:~what is it about the Nilpotence Theorem that is essentially algebro-geometric or category theoretic? In other words, given an arbitrary symmetric monoidal $\infty$-category $\C$, how many of Ravenel's conjectures can we even state, much less prove, in $\C$? Going back to Ravenel's so-called ``Orange Book" \cite{rav2}, I found that the proofs of the Nilpotence and Periodicity conjectures relied heavily on a certain sequence of spectra, called $X(n)$ by Ravenel, whose colimit was the complex cobordism spectrum $MU$. In particular, Devinatz et.~al.~showed that for any ring spectrum $R$, if an element $\alpha\in R_\ast$ is trivial in $X(n)_\ast(R)$, then it is nilpotent in $X(n-1)_\ast(R)$. From this they proved the Nilpotence Theorem, which is central to proving the Periodicity Theorem and the Thick Subcategory Theorem. 

My quest then became to understand what was really happening between $X(n)$ and $X(n+1)$. Crucial to the story here is the fact that a map of ring spectra $X(n)\to R$ determines a formal group law structure on $R_\ast$ \emph{modulo degree $n$ terms}. Thus they interpolate between the unit map $\sph\to R$ and a complex orientation $MU\to R$. Andrew Salch pointed out that Lazard, in proving that the Lazard ring is polynomial \cite{laz}, had passed from the classifying object for $n$-truncated formal group laws to $n+1$-truncated formal group laws by doing a sort of deformation theory.  Jack Morava pointed out to me that this might be a Galois theoretic or descent theoretic type of construction, and suggested I look into \cite{hess} and \cite{rog}. 

I hypothesized that the unit maps $\sph\to X(n)$ were so-called Hopf-Galois extensions like $MU$ (as Rognes had noticed in \cite{rog}). It's pretty clear that this is true, since the salient feature of being a Hopf-Galois extension is admitting an equivalence $X(n)\wedge X(n)\simeq X(n)\wedge Z$ for some ``spectral bialgebra" $Z$. Since $X(n)$ is a Thom spectrum this is pretty clearly true, with $Z=\Omega SU(n)$. However, there was a problem, because Rognes had worked only with $\E_\infty$-ring spectra and the $X(n)$ were only known to be $\E_2$. It seemed like a shot in the dark, but rather than try to recover everything Rognes had done, but for $\E_n$-ring spectra, I went to Google and typed in ``Hopf-Galois extensions of associative ring spectra." To my indescribable delight, I found that Fridolin Roth, a student of Birgit Richter, had written his thesis \cite{froth} on exactly this topic!

It followed immediately from \cite{froth} that the maps $\sph\to X(n)$ were Hopf-Galois extensions of (at least $\E_1$) ring spectra. A very natural question to ask, then, is whether or not each of the intermediate maps $X(n-1)\to X(n)$ are also some kind of Hopf-Galois extensions. This turns out to be true, and the Hopf-algebras controlling the extensions are $\sph[\Omega S^{2n+1}]$. This is proven in an unpublished paper I wrote. I was able to prove this by simply calculating the homology of everything in sight, e.g.~$X(n)\wedge_{X(n-1)}X(n)$, $X(n)\wedge \sph[\Omega S^{2n+1}]$, and the maps in between them.  One of the reasons that it's unpublished is that the proof felt unsatisfactory and like it missed some very natural structure. In particular, the equivalence $X(n)\wedge_{X(n-1)} X(n)\simeq X(n)\wedge \sph[\Omega S^{2n+1}]$ looks a lot like the generators of $H_\ast(X(n-1);\ints)$ are ``canceling out" corresponding generators in $H_\ast(X(n);\ints)$ and the only thing left is the top dimensional class. Ultimately, \cite{beardsrelthom} was an attempt to see this structure on the level of the spectra themselves, rather than just calculationally. 

In my attempts to prove something like the main theorem, I spent a lot of time looking at \cite{mahringthom} and \cite{abghr}, trying to understand the Thom diagonal as a sort of ``shear map," similar to the classical group isomorphism ${G\times_{H} G\cong G\times G/H}$ for $H\vartriangleleft G$. Eventually I realized that I would get this map for free by simply producing $X(n)$ as a Thom spectrum \emph{over} $X(n-1)$. In other words, I needed to produce a map (preferably of $\E_1$-spaces) \linebreak ${\Omega S^{2n+1}\to BGL_1(X(n-1))}$ such that the colimit of the composition $\Omega S^{2n+1}\to BGL_1(X(n-1))\to LMod_{X(n-1)}$ is equivalent to $X(n)$. Because the theory of \cite{abghr} produces a Thom isomorphism and Thom diagonal for Thom spectra over \emph{any} ring spectrum, this would give me the structure I wanted. 

Rereading (parts of) \cite{abghr}, I started to really understand the way in which they were describing Thom spectra as \emph{quotients}. It finally clicked that we can think of $X(n-1)$ as $\sph/\Omega^2 SU(n-1)$, and $X(n)$ as $\sph/\Omega^2 SU(n)$ and that the fibration $\Omega^2 SU(n)\to \Omega^2 SU(n)\to \Omega^2 S^{2n+1}$ is telling us that we should be able to go from $X(n-1)$ to $X(n)$ by just quotienting by \emph{more}. I spent some time trying to think about \emph{iterated homotopy quotients} in the form of bar constructions before I asked about this construction on MathOverflow. And just as I suspected, someone had indeed thought about this before. In particular, Qiaochu Yuan pointed out to me that these quotients can be described as Kan extensions and that iterated quotients, i.e. $(\sph/H)/(G/H)$ for a group $G$ with a subgroup $H$, can be described as iterated Kan extensions \cite{qcMO}. I had to use Lurie's \emph{operadic} Kan extensions to retain the $\E_n$-monoidal structure, but from that point on it was mostly just working out technical quasicategorical details!


\section{Colloquial summary} \label{section-summary}

In the long run, the specific result of \cite{beardsrelthom} is nice, but it's not the thing that gets me excited about algebraic topology. Rather, it's a small piece of what I consider to be a beautiful story. I described it a little in earlier sections; it's the story of \emph{chromatic homotopy theory.} I won't go into detail about it again, but the main idea is this: when we start to study things like spheres and shapes in really high dimensions, we start to see patterns arise that are coming from number theory. In other words, structures that show up when one deeply studies prime numbers \emph{also} show up when one deeply studies blobs in space.

To me, this is the really beautiful thing about mathematics: mathematicians are blindly groping in the dark, searching for structure, and every once in a while they find the exact same beautiful structure in two different places at once. It'd be like telling one engineer to build a plane, and another engineer to build a superconductor and having them separately build the exact same thing, i.e.~an object that is simultaneously a plane and a superconductor. (I don't actually know what a superconductor is, so sorry if this metaphor doesn't make any sense). Of course, we know this could never happen, since these two objects serve such different purposes, and what's more, we understand their purposes relatively well. 

However, when it comes to understanding math, we barely have any idea what we're doing. We just plow forward, mostly blind, and sometimes with a vague road map. Sometimes, like in chromatic homotopy theory, we can find connections between two \textit{a priori} disconnected things, but we still struggle to come close to understanding \emph{why} this is happening. Understanding \emph{why} would require not just knowing that I and another blind mathematician happened upon the same spot at the same time, but being able to zoom out and see how that spot relates to all the other spots. And usually we can't do this. 

This paper, and much of my work in general, is an expression of my desire to zoom out and see the whole structure. To extend our metaphor a little further, I want to trace the path of each mole and show that secretly, they were taking the same path all along. In other words, the really beautiful theorems to me are the ones that say ``Of course you're seeing structure $S$ in topic $A$ and topic $B$, because topic $B$ is actually just a sub-topic of $A$ and we never knew it!" I am of course biased, so my ultimate goal would be to say that the number theoretic structure we're seeing in chromatic homotopy theory is arising because number theory (in particular arithmetic and algebraic geometry) is itself actually a subfield of algebraic topology. Further, I think many mathematicians would like to show that \emph{many} topics of mathematics, e.g. combinatorics, algebraic geometry, topology, are all secretly different facets of the same crystal. 

The reader familiar with physics will probably have heard the term ``grand unified theory" before. The idea is that there are theories describing different physical phenomena, like quantum mechanics, and relativity, but those theories are not mutually compatible. Physicists would like for there to be a single framework that encompasses quantum mechanics, classical mechanics, relativity, and every other physical theory we know of. The ultimate goal of my mathematical research can be thought of as the same type of thing. I would love to find a since source from which the similar structures of algebraic topology and algebraic number theory both spring.




\begin{bibdiv}
\begin{biblist}

\end{biblist}
\end{bibdiv}


\begin{bibdiv}
\begin{biblist}

\bib{abg}{misc}{
Author = {Ando, Matthew}, Author= {Blumberg, Andrew J.}, Author={Gepner, David},
Title = {Parametrized spectra, multiplicative {T}hom spectra, and the twisted {U}mkehr map},
Year = {2015},
note = {arxiv.org/abs/1112.2203},
}

\bib{abghr}{article}{
    AUTHOR = {Ando, Matthew}, Author={Blumberg, Andrew J.}, Author={Gepner, David},
             Author={Hopkins, Michael J.} Author={Rezk, Charles},
     TITLE = {An {$\infty$}-categorical approach to {$R$}-line bundles,
              {$R$}-module {T}hom spectra, and twisted {$R$}-homology},
   JOURNAL = {J. Topol.},
  FJOURNAL = {Journal of Topology},
    VOLUME = {7},
      YEAR = {2014},
    NUMBER = {3},
     PAGES = {869--893}, 
}

\bib{beardsrelthom}{article}{
AUTHOR = {Beardsley, Jonathan},
TITLE = {Relative Thom Spectra Via Operadic Kan Extensions},
JOURNAL = {Algebraic and Geometric Topology},
YEAR = {2017},
Volume = {17},
NUMBER = {2},
PAGES = {1151-1162}
}

\bib{blumthh}{misc}{
Author = {Andrew J. Blumberg},
Title = {{THH} of {T}hom spectra that are $E_\infty$ ring spectra},
Year = {2008},
note = {arXiv:0811.0803},
} 

\bib{boardvogt}{book}{
AUTHOR = {Boardman, Michael}, Author={Vogt, Rainer},
TITLE = {Homotopy Invariant Algebraic Structures on Topological Spaces},
SERIES = {Lecture Notes in Mathematics},
VOLUME = {347},
PUBLISHER = {Springer-Verlag},
YEAR = {1973},
}




\bib{brown}{article}{
    AUTHOR = {Brown Jr., Edgar H.},
     TITLE = {Twisted tensor products {I}},
  JOURNAL = {Annals of Mathematics. Second Series},
    VOLUME = {69},
      YEAR = {1959},
     PAGES = {223--246},
     SERIES = {2}
}


\bib{devhopsmith}{article}{
AUTHOR = {Devinatz, E. S.}, Author={Hopkins, M. J.}, Author={Smith, Jeff},
TITLE = {Nilpotence and stable homotopy theory {I}},
JOURNAL = {Annals of Mathematics},
SERIES = {2},
VOLUME = {128},
YEAR = {1988},
PAGES = {207--241},
NUMBER = {2},
}

\bib{hess}{article}{
	AUTHOR = {Hess, Kathryn},
	TITLE = {A General Framework for Homotopic Descent and Codescent},
	eprint = {arxiv.org/abs/1001.1556v3},
	YEAR = {2010},
}


\bib{hopsmith}{article}{
	AUTHOR = {Hopkins, Michael J. }, AUTHOR = {Smith, Jeffrey H. },
	TITLE = {Nilpotence and stable homotopy theory {II}},
	JOURNAL = {Ann. of Math. (2)},
	FJOURNAL = {Annals of Mathematics. Second Series},
	VOLUME = {148},
	YEAR = {1998},
	NUMBER = {1},
	PAGES = {1--49},
}



\bib{hp}{article}{
	AUTHOR = {Hovey, Mark}, AUTHOR = {Palmieri, John H.},
	TITLE = {The structure of the {B}ousfield lattice},
	BOOKTITLE = {Homotopy invariant algebraic structures ({B}altimore, {MD},
		1998)},
	SERIES = {Contemp. Math.},
	VOLUME = {239},
	PAGES = {175--196},
	PUBLISHER = {Amer. Math. Soc., Providence, RI},
	YEAR = {1999},
}

\bib{hps}{article}{
	AUTHOR = {Hovey, Mark}, AUTHOR = {Palmieri, John H.}, AUTHOR = {Strickland, Neil P.},
	TITLE = {Axiomatic stable homotopy theory},
	JOURNAL = {Mem. Amer. Math. Soc.},
	FJOURNAL = {Memoirs of the American Mathematical Society},
	VOLUME = {128},
	YEAR = {1997},
	NUMBER = {610},
}


\bib{laz}{book}{
	AUTHOR = {Lazard, Michel},
	TITLE = {Commutative formal groups},
	SERIES = {Lecture Notes in Mathematics, Vol. 443},
	PUBLISHER = {Springer-Verlag, Berlin-New York},
	YEAR = {1975},
}


\bib{ha}{misc}{
AUTHOR = {Lurie, Jacob},
TITLE = {Higher Algebra},
YEAR = {2014},
NOTE = {math.harvard.edu/\textasciitilde lurie/papers/higheralgebra.pdf},
}

\bib{htt}{book}{
    AUTHOR = {Lurie, Jacob},
     TITLE = {Higher topos theory},
    SERIES = {Annals of Mathematics Studies},
    VOLUME = {170},
 PUBLISHER = {Princeton University Press},
   ADDRESS = {Princeton, NJ},
      YEAR = {2009},
}

\bib{mahringthom}{article}{
	AUTHOR = {Mahowald, Mark},
	TITLE = {Ring spectra which are {T}hom complexes},
	JOURNAL = {Duke Math. J.},
	FJOURNAL = {Duke Mathematical Journal},
	VOLUME = {46},
	YEAR = {1979},
	NUMBER = {3},
	PAGES = {549--559},
}


\bib{rav}{book}{
AUTHOR = {Ravenel, Doug},
TITLE = {Complex cobordism and the homotopy groups of spheres},
PUBLISHER = {Academic Press},
YEAR = {1986}
}

\bib{rav2}{book}{
	title={Nilpotence and Periodicity in Stable Homotopy Theory},
	author={Ravenel, Doug},
	series={Annals of mathematics studies},
	year={1992},
	publisher={Princeton University Press}
}


\bib{rav3}{article}{
	AUTHOR = {Ravenel, Doug},
	TITLE = {Localization with respect to certain periodic homology
		theories},
	JOURNAL = {Amer. J. Math.},
	FJOURNAL = {American Journal of Mathematics},
	VOLUME = {106},
	YEAR = {1984},
	NUMBER = {2},
	PAGES = {351--414},
}

\bib{rog}{article}{
    AUTHOR = {Rognes, John},
     TITLE = {{Galois} extensions of structured ring spectra. {S}tably
              dualizable groups},
   JOURNAL = {Mem. Amer. Math. Soc.},
  FJOURNAL = {Memoirs of the American Mathematical Society},
    VOLUME = {192},
      YEAR = {2008},
    NUMBER = {898},
}

\bib{froth}{thesis}{
	AUTHOR = {Roth, Fridolin},
	TITLE = {{Galois and Hopf-Galois} Theory for Associative {S}-Algebras},
	SCHOOL = {Universit\"{a}t Hamburg},
	YEAR = {2009},
}


\bib{qcMO}{misc}{
	AUTHOR = {Yuan, Qiaochu},
	NOTE = {Answer to MathOverflow question \textbf{218866}: \emph{Iterated Homotopy Quotient}},
	YEAR = {2015}
	}



\end{biblist}
\end{bibdiv}
\end{document}